\documentclass[reqno,12pt]{amsart}
\usepackage[sort&compress,numbers]{natbib}% bibliografia ayuda a  colocar [5-20] y que no parezca [5,6,...,20]
\usepackage{amsfonts}
\usepackage{mathdots,amssymb,graphicx,amssymb,amsmath,amsopn, mathtools,mathabx}
\usepackage[utf8]{inputenc}
\usepackage{enumitem}
\usepackage{tikz-cd}
\usepackage{scalerel}
\usepackage{pgfplots}
\usepackage{amscd}
\usepackage{latexsym}
\usepackage{xcolor}
\usepackage{amscd}
\usepackage{pb-diagram}
\usepackage{mathrsfs}
\usepackage{accents}
\usepackage[textwidth=15cm,textheight=22.5cm,headheight=0.6cm,headsep=1cm,centering]{geometry}
\usepackage{fancyhdr}
\usepackage{float}

\usepackage[bookmarks=true,colorlinks=true,linkcolor=airforceblue,citecolor=brightmaroon,plainpages=tfalse,pdftex]{hyperref}
\hypersetup{bookmarksopen=true,pdfview=fitbh}

\definecolor{brightmaroon}{rgb}{0.76, 0.13, 0.28}
\definecolor{airforceblue}{rgb}{0, 0.25, 0.77}
\definecolor{myOrange}{rgb}{1,0.5,0}
\definecolor{brightmaroon}{rgb}{0.76, 0.13, 0.28}
\definecolor{airforceblue}{rgb}{0, 0.4, 0.66}

\allowdisplaybreaks
\pgfplotsset{compat=1.18}

\theoremstyle{plain}

\newtheorem{teo}{Theorem}[section]

\newtheorem{pro}[teo]{Proposition}
\newtheorem{coro}[teo]{Corollary}
\newtheorem{defi}[teo]{Definition}
\newtheorem{remark}[teo]{Remark}

\numberwithin{equation}{section}
\newcommand{\Su}{\mathbf{S}}
\newcommand{\un}{\mathbf{u}}

\newcommand{\An}{\mathcal{A}}
\newcommand{\Bn}{\mathcal{B}}
\newcommand{\N}{\mathbb{N}}
\newcommand{\Pn}{\mathbf{P}}

\numberwithin{equation}{section}
\newcommand{\prodint}[1]{\left\langle{#1}\right\rangle}

\usepackage{color}

\newcommand*\sn{\ensuremath{\boldsymbol s}}
\newcommand*\pFq[2]{{}_{#1}F_{#2}}

\begin{document}
	
\title[Generalized Gauss-Rys orthogonal polynomials]{Generalized Gauss-Rys orthogonal polynomials}
	
\author[J. C García-Ardila,  Marcellán]{Juan C. García-Ardila,  Francisco Marcellán}
\address[J. C. Garc\'ia-Ardila]{Departamento de Matem\'atica Aplicada a la Ingenier\'ia Industrial \\Universidad Polit\'ecnica de Madrid\\ Calle Jos\'e Gutierrez Abascal 2, 28006 Madrid, Spain.}\email{juancarlos.garciaa@upm.es}
\address[F. Marcellán]{Departamento de Matemáticas, Universidad Carlos III de Madrid, Leganés, Spain} \email{pacomarc@ing.uc3m.es}
	
\thanks{ }
\date{\today}
\begin{abstract}
 Let $(P_n(x;z;\lambda))_{n\geq 0}$ be the sequence of monic orthogonal polynomials with respect to the symmetric linear functional $\sn$ defined by
$$\prodint{\sn,p}=\int_{-1}^1 p(x)(1-x^2)^{(\lambda-1/2)} e^{-zx^2}dx,\qquad\lambda>-1/2, \quad z>0.$$
In this contribution, several properties of the polynomials $P_n(x;z;\lambda)$  are studied taking into account the relation between the parameters of the three-term recurrence relation that they satisfy. Asymptotic expansions of these coefficients are given.  Discrete Painlev\'e and Painlev\'e equations associated with such coefficients appear naturally. An electrostatic interpretation of the zeros of such polynomials as well as the dynamics of the zeros in terms of the parameters $z$ and $\lambda$ are given.
\end{abstract}
	
\maketitle
\tableofcontents
\section{Introduction}
	In this contribution, we deal with sequences of monic polynomials that are orthogonal with respect to the weight function $\omega(x; z)= (1-x^2)^{\lambda-1/2} e^{-z x^{2}},$ $ x\in (-1,1)$ with $\lambda >-1/2$, $z\geq0,$  which is a product of the Gegenbauer weight function and the Rys exponential function. These polynomials have been studied in \cite{Gradimir22} in the framework of Gauss-quadrature rules. Notice that when $z=0$ you get the Gegenbauer weight and when $\lambda=1/2$ the so-called Rys weight appears. The last case has been extensively studied in the literature. More precisely,  Rys quadrature formulas were introduced in \cite{Du76}  as a method in computational quantum chemistry to evaluate two-electron repulsion integrals, which appear in molecular quantum mechanical calculations involving Gaussian Cartesian basis functions. In \cite{K16,Shizgal} the existing strategies for the evaluation of Rys nodes and weights were analyzed. An efficient and stable method for constructing Gauss-Rys quadratures has recently been given in \cite{Gradimir18}.\\
	
	Here we focus the attention on the properties of the sequence of monic  polynomials $(P_{n})_{n\geq0}$ orthogonal with respect the weight function $\omega(x; z).$ They satisfy a three-term recurrence relation (TTRR)
	\begin{equation}\label{rr}
		\begin{aligned}
			x\,P_{n}(x)&=P_{n+1}(x)+ \gamma_{n}\,P_{n-1}(x),\quad n\geq0,\\
			P_{-1}(x)&=0, \ \ \ \ \ P_{0}(x)=1.
		\end{aligned}
	\end{equation}
	
	Taking into account the semiclassical character of the weight function, we are able to deduce a non linear difference equation that the parameters $\gamma_{n}$ satisfy. It is a discrete Painlev\'e equation. On the other hand, since such parameters depend on $z$, a differential Painlev'e equation is deduced. This result is strongly related with the fact that you have a continuous perturbation of the Gegenbauer weight that yields a Volterra lattice, also known in the literature as the discrete KdV equation and Langmuir lattice. It is a system of ordinary differential equations with variables indexed by some of the points of a $1$-dimensional lattice. The Volterra lattice is a special case of the generalized Lotka–Volterra equation describing predator–prey interactions for a sequence of species such that each specie is preying on the next one in the sequence.  The Volterra lattice is an integrable system related to the Toda lattice. It is also used as a model for Langmuir waves in plasma.\\
	
	Notice that when $\lambda=1/2$ these polynomial sequences are related to truncated Hermite orthogonal polynomials, i. e., the corresponding measure is the truncated normal distribution in a symmetric interval of the real line. The analytic properties of such polynomials as well as the electrostatic interpretation of their zeros have been studied in \cite{DM22}. Truncated generalized Hermite polynomials have been studied in \cite{DGM23} by using a symmetrization process from the sequence of orthogonal polynomials associated with the truncated gamma distribution (see also \cite{GM23}). Notice that \cite{Chen05} is a seminal paper on orthogonal polynomials with respect to discontinuous weights. The techniques used therein are very useful in the analysis of Gaussian unitary ensembles with jump discontinuities in random matrices \cite{Lyu1}.\\
	
	The framework of our approach is the theory of semiclassical linear functionals (for more information, see \cite{Ma87,Ma91}). Their connections with discrete and differential equations of Painlev\'e type satisfied by the coefficients of the three-term recurrence relation of the corresponding sequences of orthogonal polynomials have attracted the interest of many researchers see \cite{Magnus1,Magnus2,Walter}, among others.\\
	
	To construct the Gaussian formula with at most $N$ nodes, using the classical Golub-Welsch algorithm, we need the first $N$  coefficients $\gamma_{n}$ in  \eqref{rr},  $n=0,1\ldots,N-1,$ which must be computed numerically for this semiclassical case. Such approaches belong to the so-called constructive theory of orthogonal polynomials developed in \cite{Ga82,Ga04,Gradimir14}. In the numerical construction of the recursive coefficients, an important aspect is the sensitivity of the problem concerning a small perturbation in the input. However, recent progress in variable-precision arithmetic and symbolic computation enables us to generate recursive coefficients directly by applying the original Chebyshev method of moments but by using sufficiently high-precision arithmetic to overcome the numerical instability.\\
	
	The structure of this paper is as follows. In Section~\ref{section2}, a basic background about orthogonal polynomials and linear functionals is given. We emphasize the properties of semiclassical linear functionals, which will play a central role in the sequel. Section~\ref{section3} is focused on the fact that the measure associated with generalized Gauss-Rys orthogonal polynomials is semiclassical of class 2. As a consequence, an explicit expression of the moments in terms of an $_{1}F_{1}$ confluent hypergeometric functions is given. On the other hand, a first order linear differential equation that the Stieltjes function associated with the measure satisfies is given. In Section~\ref{section4}, the Laguerre-Freud equation satisfied by the parameters $\gamma_{n}$ of the TTRR is obtained. As a consequence, after a change in the above parameters, a modified discrete Painlev\'e~II equation satisfied by the new parameters is deduced. Section \ref{section5} deals with the analysis of ladder operators for generalized Gauss-Rys orthogonal polynomials. Thus, a second order linear differential equation they satisfy follows in a natural way. In Section~\ref{Section6}, the behavior of the parameters $\gamma_{n}$ in terms of $z$ is given.  A second order nonlinear differential equation of Painlev\'e type is given. Finally, in Section~\ref{section7} from the holonomic equation we get an electrostatic interpretation of the zeros of generalized Gauss-Rys polynomials. Its dynamic behavior is also pointed out.
	%%%%%%%%%%%%%%%%%%%%%%%%%%%%%%%%%%%%%%%
	%%%%%%%%%%%%%%%%%%%%%%%%%%%%%%%%%%%%%%%
	\section{Basic background}\label{section2}
	Let $\un$ be a linear functional defined on the linear space of polynomials with complex coefficients $\mathbb{P}$, that is, $$\un:\mathbb{P}\to \mathbb{C}, \quad p(x)\to \prodint{\un,p(x)}.$$ Denote the $n$th moment of $\un$ by $\un_n:=\prodint{\un,x^n}$, $n\geq0$.
	\begin{defi} Let $\un$  be a linear functional.
		\begin{enumerate}
			\item For $q\in\mathbb{P}$ we define the linear functional $q(x)\un$ as 	\begin{equation*}
				\prodint{q(x)\un,p(x)}=\prodint{\un,q(x)p(x)}, \quad p\in\mathbb{P}.
			\end{equation*}
			\item The derivative of  $\un$ is the linear functional $D\un$ such that
			$$\prodint{D\un,p(x)}=-\prodint{\un,p^\prime(x)}.$$
		\end{enumerate}
	\end{defi}
	The linear functional $\un$ is said to be quasi-definite (resp. positive-definite) if every leading principal submatrix of the Hankel matrix $H=(\un_{i+j})_{i , j=0}^{\infty}$ is nonsingular (resp. positive-definite).  In such a situation, there exists a sequence of monic polynomials $(P_n)_{n\geq0}$ such that $\deg{P_n}=n$ and $\prodint{\un,P_n(x)P_m(x)}=K_{n}\delta_{n,m},$ where $\delta_{n,m}$ is the Kronecker symbol and $K_{n}\ne 0$ (see \cite{Chi,GMM21}). The sequence $(P_n)_{n\geq0}$ is said to be the sequence of monic orthogonal polynomials (SMOP) with respect to the linear functional~$\un$. \\
	
	Let $\un$ be a quasi-definite linear functional and let $(P_n)_{n\geq0}$ be its corresponding SMOP. Then  there  exist two  sequences of complex numbers $( a_n)_{n\geq1}$ and $(b_n)_{n\geq0}$, with $ a_n\ne 0$, such that	
	\begin{equation}\label{ttrrr}
		\begin{aligned}
			x\,P_{n}(x)&=P_{n+1}(x)+b_n\,P_{n}(x)+ a_{n}\,P_{n-1}(x),\quad n\geq0,\\
			P_{-1}(x)&=0, \ \ \ \ \ P_{0}(x)=1.
		\end{aligned}
	\end{equation}
	Conversely, by Favard's Theorem \cite{Chi,GMM21}, if $(P_n)_{n\geq0}$ is a sequence of monic polynomials generated by a three-term recurrence relation (TTRR) as in \eqref{ttrrr} with $ a_n\ne 0,$ $n\geq1$, then there exists a unique linear functional $\un$ such that $(P_n)_{n\geq0}$ is its corresponding SMOP.
	\begin{defi}
		Given a quasi-definite linear functional $\un$ with moments $(\un_n)_{n\geq0}$, the formal series
		\begin{equation*}
			\Su_{\un}(z)=:\sum_{n=0}^\infty \dfrac{\un_n}{z^{n+1}}
		\end{equation*}
		is said to be the \textit{Stieltjes function} associated with $\un$.
	\end{defi} Due to the connection between the Stieltjes function, orthogonal polynomials (as denominators of the diagonal  Pad\'e approximants to such functions)  and the coefficients of the three-term recurrence relation \eqref{ttrrr} through   continued fractions, the Stieltjes function is a very useful tool to study properties of quasi-definite linear functionals \cite{MR1164865,MR2442472,MR1083352,MR0025596}.\\
	
	Next we summarize some basic results concerning semiclassical linear functionals, which will play a central role in the sequel.
	\begin{defi}[\cite{Ma87}]
		A quasi-definite functional $\mathbf{u}$ is said to be semiclassical if there exist non-zero polynomials $\phi(x)$ and $\psi(x)$ with $\deg\phi(x)=:r\ge 0$ and $\deg\psi(x)=:t\ge 1$, such that $\mathbf{u}$ satisfies the distributional Pearson equation
		\begin{equation}\label{pearson-semic}
			D(\phi(x)\,\mathbf{u})+\psi(x)\,\mathbf{u}=0.\end{equation}
		A sequence of orthogonal polynomials associated with $\mathbf{u}$ is called a semiclassical sequence of orthogonal polynomials.
	\end{defi}
	Notice that a semiclassical linear functional satisfies many Pearson equations, Indeed, if $q(x)$ is a polynomial and $\un$ satisfies \eqref{pearson-semic}, then $\un$ also satisfies the Pearson equation
	$$D(\widetilde\phi(x)\,\mathbf{u})+\widetilde\psi(x)\,\mathbf{u}=0,$$
	where $\widetilde\phi(x)=q(x)\phi(x)$ and $\widetilde\psi(x)=(\psi(x)q(x)-\phi(x)q^\prime(x)).$
	The non-uniqueness of the Pearson equation motivates the following definition.
	\begin{defi}
		The class of a semiclassical functional $\mathbf{u}$ is defined as
		\begin{equation*}%\label{class}
			\mathfrak{s}(\mathbf{u})=: \min \max\{\deg \phi(x)-2, \deg\psi(x)-1 \},
		\end{equation*}
		where the minimum is taken among all pairs of polynomials $\phi(x)$ and $\psi(x)$ so that \eqref{pearson-semic} holds.
	\end{defi}
	\begin{pro}[\cite{GMM21, Ma91}]\label{sim_cond}
		Let $\mathbf{u}$ be a semi-classical linear functional and let $\phi(x)$ and $\psi(x)$ be non-zero polynomials with $\deg\phi(x)=:r$ and $\deg \psi(x)=:t$, such that \eqref{pearson-semic} is satisfied. Let $s =: \max(r-2,t-1)$. Then $s = \mathfrak{s}(\mathbf{u})$  if and only if
		\begin{equation}\label{prodformula}\prod_{c:\,\phi(c)=0}\left(|\psi(c)+\phi^\prime(c)|+|\langle\mathbf{u},\theta_c\psi(x)+\theta^2_c\phi(x)\rangle|\right)>0.
		\end{equation}
		Here, $\theta_c f(x)=\dfrac{f(x)-f(c)}{x-c}.$
	\end{pro}
	
	According to the class, there is a hierarchy of semiclassical linear functionals. The class zero is constituted by the classical linear functionals (Hermite, Laguerre, Jacobi and Bessel). The class one has been studied in  \cite{Belmehdi}.
	%\begin{defi}A quasi-definite linear  functional  $\mathbf{\un}$ is said to be \textit{symmetric} if
	%$$\prodint{\un,x^{2k+1}}=0,\quad k\ge 0,$$
	%that is, all its odd moments are zero.
	%\end{defi} 	
	%The following result characterizes symmetric quasi-definite linear functionals in terms of the recurrence relation satisfied by its MOPS.
	%\begin{teo}[\cite{Chi}]\label{symmetric}
	%Let $(S_n(x))_{n\ge 0}$ be an MOPS associated with a quasi-definite linear functional~$\mathbf{\un}$. Then the following statements are equivalent:
	%\begin{enumerate}
	%\item[(i)] $\mathbf{\un}$ is symmetric.
	%\item[(ii)] $S_{n}(x)$ has the same parity as $n$, that is, $S_{n}(-x)= (-1)^{n} S_{n}(x).$
	%\item[(iii)] $(S_n(x))_{n\ge 0}$ satisfies a recurrence relation
	%\begin{equation*}%\label{Tre}
	%\begin{aligned}
	%xS_n(x)&=S_{n+1}(x)+\gamma_n\,S_{n-1}(x), \ \ \ \ \ n\geq0,\\
	%S_{-1}(x)&=0, \ \ \ \ \ S_{0}(x)=1.
	%\end{aligned}
	%\end{equation*}
	%\end{enumerate}
	%	\end{teo}
%Let $S_{2n}(x)=P_{n}(x^2), S_{2n+1}(x)= xQ_{n}(x^2). $
%\begin{pro}[Chihara \cite{Chi}]%\label{corosymm}
%Let $\mathbf{s}$  be a quasi-definite  symmetric linear functional. Then the sequence of monic polynomials  $(P_{n}(x))_{n\ge 0}$  is orthogonal with %respect to the functional  $\un$  defined by  $\prodint{\un,x^n}=\prodint{\sn,x^{2n}}$ and the sequence of monic polynomials$(Q_{n}(x))_{n\ge 0}$ is %orthogonal with respect to the functional $x\un$.
%	\end{pro}
\section{Generalized Gauss-Rys orthogonal polynomials}\label{section3}
In this section we study the Generalized Gauss-Rys linear functional $\sn$,
\begin{equation}\label{innerproduct}
\prodint{\sn,p}=\int_{-1}^1p(x)(1-x^2)^{(\lambda-1/2)} e^{-zx^2}dx,\qquad\lambda>-1/2, \quad z>0.
\end{equation}
From a computational point of view, in order to analyze quadrature rules for such weight functions, zeros and weights have been recently analyzed in \cite{Gradimir22}.
We will prove that $\sn$ is a semiclassical functional of class $s=2$.	Indeed, we will first find the corresponding Pearson equation. In the sequel, we will deduce a second-order linear recurrence equation for its moments. We will also get a first-order linear differential equation satisfied by the Stieltjes function.
\begin{teo}
The functional $\sn$ is semi-classical of class 2. Moreover, $\sn$ satisfies the Pearson equation \begin{equation}\label{pearson_eq}
	D(\phi\sn)+\psi\sn=0,
\end{equation}
with $\phi(x)=(1-x^2)$ and $\psi(x)=(2z+2\lambda+1)x-2zx^3.$
\end{teo}
\begin{proof} Using integration by parts
\begin{align*}
	\prodint{D(\phi(x)\sn),p(x)}&=-\prodint{\sn,\phi(x) p^\prime(x)}=-\int_{-1}^1\phi(x) p^\prime(x)(1-x^2)^{\lambda-1/2} e^{-zx^2}dx\\
	&=-p(x)(1-x^2)^{\lambda+1/2}e^{-zx^2}\left.\right|_{-1}^{1}\\
	&+\int_{-1}^1 p(x)\left(- 2 x -2(\lambda-1/2)x -2zx (1-x^2)\right)(1-x^2)^{\lambda-1/2} e^{-zx^2}dx \\
	&=-\prodint{\sn,\psi(x)p(x)}.
\end{align*}
\end{proof}
\begin{pro}
Let $(\sn_m)_{m\geq 0}$ be the sequence of moments of $\sn$. Then
\begin{equation}\label{moments}
	\sn_m=\begin{cases}
		\dfrac{
			\Gamma(n+1/2)\Gamma(\lambda+1/2)}{\Gamma(\lambda+n+1)}\times\pFq{1}{1}\left(n+1/2;\lambda+n+1,-z\right),& m=2n,\\
		0,&m=2n+1,
	\end{cases}\quad n\geq 0,
\end{equation}
where
$$\pFq{p}{q}\left(a_1,\ldots,a_p;b_1,\ldots,b_q,x\right)=\sum_{k=0}^\infty\dfrac{(a_1)_k\cdots(a_p)_k}{(b_1)_k\cdots(b_q)_k}\dfrac{x^k}{k!} $$ and
$(a)_k$ is the \textit{Pochhammer symbol} defined by $(a)_0=1,$
\begin{equation}\label{poh}
	(a)_k=a(a+1)\cdots (a+k-1)=\dfrac{\Gamma(a+k)}{\Gamma(a)},	\quad k=1,2\ldots
\end{equation}
\end{pro}
\begin{proof}
Since $(1-x^2)^{\lambda-1/2} e^{-zx^2}$ is a symmetric weight, then $$\sn_{2n+1}=\prodint{\sn,x^{2n+1}}=0$$ for all $n\geq 0.$
On the other hand
\begin{equation}\label{33}
	\begin{aligned}
		\sn_{2n}=\prodint{\sn,x^{2n}}&=2\int_0^1x^{2n}(1-x^2)^{\lambda-1/2} e^{-zx^2}dx\\
		&=\int_0^1u^{n-1/2}(1-u)^{\lambda-1/2} e^{-zu}du.	
	\end{aligned}
\end{equation}
Thus, using the integral formula for the Kummer confluent hypergeometric function \cite[pag. 326]{OLB10}
$$\pFq{1}{1}\left(a;b,z\right)=\dfrac{\Gamma(b)}{\Gamma(a)\Gamma(b-a)}\int_0^1x^{a-1}(1-x)^{b-a-1}e^{zx}dx,\quad \operatorname{Re}(b)>\operatorname{Re}(a)>0,$$
and taking $a=n+1/2$ and $b=\lambda+n+1$ in \eqref{33} we get the result.
\end{proof}
\begin{pro}
Let $(\sn_n)_{n\geq 0}$ be the sequence of moments of $\sn$. They satisfy the recurrence relation
\begin{equation}\label{recmoments}
	2z\sn_{n+3}-(n+2z+2\lambda+1)\sn_{n+1}+n\sn_{n-1}=0, \quad n\geq 0,
\end{equation}
with initial conditions
\begin{equation}\label{initialcondi}
	\begin{aligned}
		\sn_0&=			\sqrt{\pi}\dfrac{\Gamma(\lambda+1/2)}{\Gamma(\lambda+1)}\times\pFq{1}{1}\left(1/2;\lambda+1,-z\right),\\[10pt] \sn_2&=			\dfrac{\sqrt{\pi}}{2}\dfrac{\Gamma(\lambda+1/2)}{\Gamma(\lambda+2)}\times\pFq{1}{1}\left(3/2;\lambda+2,-z\right).	
	\end{aligned}		
\end{equation}
\end{pro}
\begin{proof}
From the Pearson equation \eqref{pearson_eq} we have
\begin{align*}
	\prodint{D(\phi\sn)+\psi\sn,x^n}&=-\prodint{\sn,n\phi(x)\,x^{n-1}}+\prodint{\sn,\psi(x)\,x^n}\\
	&=-\prodint{\sn,n(x^{n-1}-x^{n+1})}+\prodint{\sn,2(\lambda+1/2+z)x^{n+1}-2zx^{n+3}}\\
	&=n\sn_{n+1}-n\sn_{n-1}+2(\lambda+1/2+z)\sn_{n+1}-2z\sn_{n+3}.
\end{align*}
To deduce the class of the functional, we will use Proposition \ref{sim_cond}. In our case \eqref{prodformula} becomes
$$\prod_{c:\,\phi(c)=0}\left(\left|2\lambda-1\right|+\left|(2\lambda+2z-2cz)\sn_0-2z\sn_2\right|\right),$$
where $c=\pm1$  are the zeros of $\phi$. Note that the above expression is different from zero if $\lambda\ne1/2$. When $\lambda=1/2$ we use the error function ${\rm erf}(z)$ defined by \cite[7.2.1]{OLB10}
$${\rm erf}(z)=\dfrac{2}{\sqrt{\pi}}\int_0^ze^{-x^2}dx,\ \text{which satisfies}\ \int{\rm erf }(z)dz=z\,{\rm erf}(z)+\dfrac{e^{-z^2}}{\sqrt{\pi}}+C,$$
where $C$ is the integration constant. Thus,
$$\begin{aligned}
	(1+2z-2cz)\sn_0-2z\sn_2=&2\,{\dfrac{z (1-c)\, {\rm erf} \left(\sqrt {z}\right)
			\sqrt {\pi}+\sqrt {z}{{e}^{-z}}}{\sqrt {z}}}
\end{aligned}$$
which is always different from zero  for $z>0.$
\end{proof}
\begin{pro}\label{Stielt}
The Stieltjes function $\Su_{\sn}(t;z;\lambda)$ associated with the linear functional $\sn$ satisfies the non-homogeneous  first order linear ODE
$$\phi(z)\partial_t\Su_{\sn}(t;z;\lambda)+[\phi^\prime(z)+\psi(z)]\Su_{\sn}(t;z;\lambda)=-2z\sn_2+((2z+2\lambda-2)t-2zt^2)\sn_0.$$
\end{pro}
\begin{proof}
First, note that
$$\Su_{\sn}(t;z;\lambda)=\sum_{n=0}^\infty\dfrac{\sn_{2n}}{t^{2n+1}}\quad\text{and}\quad\partial_t\Su_{\sn}(t;z;\lambda)=-\sum_{n=0}^\infty(2n+1)\dfrac{\sn_{2n}}{t^{2n+2}}.$$
Now, from \eqref{recmoments}
\begin{align*}
	\Su_{\sn}(t;z;&\lambda)=\dfrac{\sn_0}{t}+\dfrac{\sn_2}{t^3}+\sum_{n=2}^\infty\dfrac{\sn_{2n}}{t^{2n+1}}\\
	&=\dfrac{\sn_0}{t}+\dfrac{\sn_2}{t^3}+\dfrac{1}{2z}\sum_{n=0}^\infty(2n+2z+2\lambda+2)\dfrac{\sn_{2n+2}}{t^{2n+5}}-\dfrac{1}{2z}\sum_{n=0}^\infty(2n+1)\dfrac{\sn_{2n}}{t^{2n+5}} \\
	&=\dfrac{\sn_0}{t}+\dfrac{\sn_2}{t^3}+\dfrac{(2z+2\lambda-1)}{2zt^2}\left(\sum_{n=0}^\infty\dfrac{\sn_{2n+2}}{t^{2n+3}}+\dfrac{\sn_0}{t}-\dfrac{\sn_0}{t}\right)\\
	&+\dfrac{1}{2zt}\left(\sum_{n=0}^\infty(2n+3)\dfrac{\sn_{2n+2}}{t^{2n+5}}+\dfrac{\sn_0}{t^2}-\dfrac{\sn_0}{t^2}\right)-\dfrac{1}{2zt^3}\sum_{n=0}^\infty(2n+3)\dfrac{\sn_{2n}}{t^{2n+3}}\\
	&=\dfrac{\sn_0}{t}+\dfrac{\sn_2}{t^3}+\dfrac{(2z+2\lambda-1)}{2zt^2}\left(\Su_{\sn}(t;z;\lambda)-\dfrac{\sn_0}{t}\right)\\
	&-\dfrac{1}{2zt}\left(\partial_t\Su_{\sn}(t;z;\lambda)-\dfrac{\sn_0}{t^2}\right)+\dfrac{1}{2zt^3}\partial_t\Su_{\sn}(t;z;\lambda).
\end{align*}
Thus, multiplying both sides of the equation by $2zt^3$ and reorganizing
$$(1-t^2)\partial_t\Su_{\sn}(t;z;\lambda)+t((2z+2\lambda-1)t-2zt^2)\Su_{\sn}(t;z;\lambda)=-2z\sn_2+(2z+2\lambda-1-2zt^2)\sn_0.$$
\end{proof}
\begin{remark}%\label{rem3,7}
The case $\lambda=1/2$ was studied by D. Dominici and F. Marcellán in~\cite{DM22}. Thus, we shall deal with the case $\lambda\ne1/2.$
\end{remark}
\section{Laguerre-Freud's equations}\label{section4}
If $\sn$ is a semiclassical linear functional satisfying \eqref{pearson_eq}, then the coefficients of the three-term recurrence relation \eqref{ttrrr} satisfy a nonlinear system of equations obtained from the relations.
\begin{align}
\label{1eq}&\prodint{\psi\sn,P_n^2}=-\prodint{D(\phi\sn),P_n^2},\\
\label{2eq}&\prodint{\psi\sn,P_{n+1}P_n}=-\prodint{D(\phi\sn),P_{n+1}P_n},
\end{align}
where $(P_n)_{n\geq 0}$ is the sequence of monic orthogonal polynomials with respect to~$\sn$. These equations are known in the literature as \textit{Laguerre-Freud equations} (see \cite{Said,Magnus1,Magnus2}).  Suppose that  $(P_n)_{n\geq 0}$ satisfies the recurrence relation
\begin{equation}\label{ttrrs}
\begin{aligned}
	&xP_{n}(x)=P_{n+1}(x)+\gamma_{n}P_{n-1}(x), \quad n\geq 0,	\\
	&P_{-1}(x)=0,\quad P_{0}(x)=1.
\end{aligned}	
\end{equation}
Notice that from the parity of $P_n(x)$
\begin{equation}\label{rev}
P_n(x)=x^n+\lambda_{n,n-2}x^{n-2}+\lambda_{n,n-4}x^{n-4}+\mathcal{O}(x^{n-6}).
\end{equation}
Then
\begin{align*}
x&P_n(x)=x^{n+1}+\lambda_{n,n-2}x^{n-1}+\lambda_{n,n-4}x^{n-3}+\mathcal{O}(x^{n-5})\\
&=x^{n+1}+(\lambda_{n+1,n-1}+\gamma_n)x^{n-1}+(\lambda_{n+1,n-3}+\lambda_{n-1,n-3}\gamma_n)x^{n-3}+\mathcal{O}(x^{n-5}).
\end{align*}
In particular,
\begin{equation}\label{cases}
\gamma_n=\lambda_{n,n-2}-\lambda_{n+1,n-1},
\end{equation}
or, equivalently,
\begin{equation*}%\label{sumgama}
\sum_{j=0}^n\gamma_{j}=-\lambda_{n+1,n-1}.	
\end{equation*}
Now, reverse \eqref{rev} and take into account \eqref{cases}
\begin{align}
\notag\partial_xP_n(x)&=nx^{n-1}+(n-2)\lambda_{n,n-2}x^{n-3}+\mathcal{O}(x^{n-5})\\
\notag&=nP_{n-1}(x)+(-n\lambda_{n-1,n-3}+(n-2)\lambda_{n,n-2})P_{n-3}(x)+\mathcal{O}(x^{n-5})\\
\label{der}&=nP_{n-1}(x)-(n\gamma_{n-1}+2\lambda_{n,n-2})P_{n-3}(x)+\mathcal{O}(x^{n-5}).
\end{align}
From the recurrence relation \eqref{ttrrs} we get
\begin{equation}\label{rrttx^2}
x^2P_n(x)=P_{n+2}(x)+(\gamma_{n+1}+\gamma_n)P_{n}(x)	+\gamma_n\gamma_{n-1}P_{n-2}(x), \quad n\geq 0.
\end{equation}
Keeping this in mind, we can deduce the following nonlinear difference equation that the coefficients $\gamma_{n}$ of the TTRR satisfy.
\begin{pro}
The coefficients of the three-term recurrence relation \eqref{ttrrs} associated with the linear functional $\sn$ satisfy the Laguerre-Freud equation
\begin{multline}\label{111}
	\dfrac{1}{2}=\Big[z+\lambda+n+1-z\gamma_{n+2}-z\gamma_{n+1}\Big]\gamma_{n+1}-\Big[z+\lambda+n-1-z\gamma_{n}-z\gamma_{n-1}\Big]\gamma_{n}, \quad n\geq1,
\end{multline}
with $\gamma_0=0.$		
\end{pro}
\begin{proof}
Let $h_n=\prodint{\sn,P_n^2}$.
Equation \eqref{1eq} is equivalent to
\begin{equation}\label{eqq1}
	\prodint{\sn,\psi P_n^2}=2\prodint{\sn,\phi P_nP^\prime_n},
\end{equation}
where $\phi(x)=(1-x^2)$ and $\psi(x)=(1+2\lambda+2z)x-2zx^3.$ Since $P_n^2(x)$ and $\phi(x)$ are even functions and $\psi(x)$ and $P_n(x)P^\prime_n(x)$ are odd functions, then by the symmetry of $\sn$ we conclude that the left-hand side and the right-hand side of \eqref{eqq1} are zero. Therefore, this equation does not provide any information. On the other hand, \eqref{2eq} is equivalent to
\begin{equation*}%\label{eqq2}
	\prodint{\sn,\psi P_nP_{n+1}}=\prodint{\sn,\phi P_nP^\prime_{n+1}}+\prodint{\sn,\phi P^\prime_{n}P_{n+1}}.   \end{equation*}
Let $a=1+2\lambda+2z$ and $b=2z.$ Then from \eqref{der} and\eqref{rrttx^2}
\begin{equation}\label{eqq2}
	\begin{aligned}
		\prodint{\sn,\psi P_nP_{n+1}}&=a\prodint{\sn,P_n(P_{n+2}+\gamma_{n+1}P_n)}-b\prodint{\sn,x^2P_n(P_{n+2}+\gamma_{n+1}P_n)}\\
		&=a\gamma_{n+1}h_n-b\prodint{\sn,\gamma_{n+2}\gamma_{n+1}P_{n}^2}-b\prodint{\sn,\gamma_{n+1}(\gamma_{n+1}+\gamma_{n})P_{n}^2}\\
		&=(a\gamma_{n+1}-b\gamma_{n+2}\gamma_{n+1}-b\gamma_{n+1}(\gamma_{n+1}+\gamma_{n}))h_n,\\
	\end{aligned}
\end{equation}
\begin{equation}\label{eqq3}
	\begin{aligned}
		\prodint{\sn,\phi P^\prime_{n}P_{n+1}}&=-\prodint{\sn,x^2P_n^\prime P_{n+1}}=-\prodint{\sn,nx^2P_{n-1}P_{n+1}}=-n\gamma_{n+1}h_n,\end{aligned}
\end{equation}
\begin{equation*}
	\begin{aligned}	&\prodint{\sn,\phi P_nP^\prime_{n+1}}=\prodint{\sn, P_nP^\prime_{n+1}}-\prodint{\sn,x^2 P_nP^\prime_{n+1}}\\
		&=(n+1)h_n-\prodint{\sn,x^2 P_n[(n+1)P_{n}-((n+1)\gamma_{n}+2\lambda_{n+1,n-1})P_{n-2}+\mathcal{O}(x^{n-4})]}\\
		&=(n+1)h_n+((n+1)\gamma_{n}+2\lambda_{n+1,n-1})h_n-(n+1)\prodint{\sn,x^2P_n^2}\\
		&=[n+1+2\lambda_{n+1,n-1}-(n+1)\gamma_{n+1}]h_n.
	\end{aligned}
\end{equation*}
From the above and \eqref{cases}, we have
\begin{equation*}
	a\gamma_{n+1}-\gamma_{n+2}\gamma_{n+1}b-\gamma_{n+1}(\gamma_{n+1}+\gamma_{n})b=+2\lambda_{n+1,n-1}-(2n+1)\gamma_{n+1}+n+1.
\end{equation*}
Shifting $n\to n-1$, subtracting and taking into account \eqref{cases} \begin{multline*}
	a(\gamma_{n+1}-\gamma_{n})+\gamma_{n+1}(\gamma_{n}-\gamma_{n+2})b+\gamma_{n}(\gamma_{n}+\gamma_{n-1})b-\gamma_{n+1}(\gamma_{n+1}+\gamma_{n})b\\=(2n-3)\gamma_{n}-(2n+1)\gamma_{n+1}+1.
\end{multline*}
According to the values of $a$ and $b$ we get \eqref{111}.
\end{proof}
\begin{remark}
Notice that the above result states a quadratic difference equation involving four consecutive terms of the sequence $(\gamma_{n})_{n\geq0}.$
\end{remark}
Our next step will be to deduce a cubic difference equation involving three consecutive terms of a new sequence defined from $(\gamma_{n})_{n\geq0}.$

Let $v(x)=zx^2-(\lambda-1/2)\ln(1-x^2)$ which is two times derivable on $(-1,1)$. Moreover,
\begin{equation*}
v^\prime(x)=2zx+(2\lambda+1)\dfrac{x}{1-x^2}.
\end{equation*}
If we define $\An_n(x)$ and $\Bn_n(x)$ by
\begin{align*}
\An_{n}(x)&=\dfrac{1}{h_n}\int^1_{-1}\dfrac{v^\prime(x)-v^\prime(y)}{(x-y)}P_n^2(y)w(y)dy,\\
\Bn_{n}(x)&=\dfrac{1}{h_{n-1}}\int^1_{-1}\dfrac{v^\prime(x)-v^\prime(y)}{(x-y)}P_n(y)P_{n-1}(y)w(y)dy,
\end{align*}
then (see \cite[Chapter 3]{Mourad} as well as \cite[Theorem 4.2]{Walter})
\begin{align}
\label{B1}\Bn_{n+1}(x) + \Bn_n(x)&= x\An_n(x) - v^\prime(x),\\
\label{B2}\gamma_{n+1}\An_{n+1}(x) - \gamma_n \An_{n-1}(x)&= 1 + x \left[\Bn_{n+1}(x)-\Bn_n(x)\right].
\end{align}
In our case
$$\dfrac{v^\prime(x)-v^\prime(y)}{(x-y)}=2z+\dfrac{2\lambda-1}{2(1-x)(1-y)}+\dfrac{2\lambda-1}{2(1+x)(1+y)}.$$
Therefore if $w(x)=(1-x^2)^{\lambda-1/2}e^{-zx^2}$, then
\begin{align*}
\An_n(x)=2z+\dfrac{R_n}{1-x}+\dfrac{r_n}{1+x},\qquad
\Bn_n(x)=\dfrac{T_n}{1-x}+\dfrac{t_n}{1+x},
\end{align*}
where
\begin{align*}
R_n&=\dfrac{2\lambda-1}{2h_n}\int_{-1}^1\dfrac{P_n^2(y)}{(1-y)}w(y)dy,  & r_n&=\dfrac{2\lambda-1}{2h_n}\int_{-1}^1\dfrac{P_n^2(y)}{(1+y)}w(y)dy, \\
T_n&=\dfrac{2\lambda-1}{2h_{n-1}}\int_{-1}^1\dfrac{P_n(y)P_{n-1}(y)}{(1-y)}w(y)dy,  & t_n&=\dfrac{2\lambda-1}{2h_{n-1}}\int_{-1}^1\dfrac{P_n(y)P_{n-1}(y)}{(1+y)}w(y)dy.
\end{align*}
Thus, from  \eqref{B2}
\begin{multline*}
\gamma_{n+1}\left[2z(1-x^2)+(1+x)R_{n+1}+(1-x)r_{n+1}\right]\\ -\gamma_{n}\left[2z(1-x^2)+(1+x)R_{n-1}+(1-x)r_{n-1}\right]\\
=1-x^2+x\left((1+x)T_{n+1}+(1-x)t_{n+1}\right)-x\left((1+x)T_{n}+(1-x)t_{n}\right).
\end{multline*}Comparing the coefficients  we get the relations
\begin{align}
&\gamma_{n+1}(2z+R_{n+1}+r_{n+1})-\gamma_{n}(2z+R_{n-1}+r_{n-1})=1,\label{eq1}\\
& \gamma_{n+1}(R_{n+1}-r_{n+1})-\gamma_{n}(R_{n-1}-r_{n-1})=(T_{n+1}-T_{n})+(t_{n+1}-t_{n}),\label{eq2}\\
&2z(\gamma_{n+1}-\gamma_{n})=1-(T_{n+1}-T_{n})+(t_{n+1}-t_{n}).  \label{eq3}
\end{align}
On the other hand, from \eqref{B1} we get
\begin{equation*}
(1+x)T_{n+1}+(1-x)t_{n+1}+(1+x)T_{n}+(1-x)t_{n}=(1+x)xR_{n}+(1-x)xr_{n}-(2\lambda-1)x   \end{equation*}
and again, comparing coefficients
\begin{align}
\notag&T_{n+1}+t_{n+1}+T_{n}+t_{n}=0 ,\\
\label{eq5}&T_{n+1}-t_{n+1}+T_{n}-t_{n}=R_n+r_n-(2\lambda-1),\\
\label{eq6}&R_n-r_n=0.
\end{align}
\begin{teo} For $n\geq 1$, the coefficients $\gamma_n(z;\lambda)$ satisfy the equation
\begin{equation}\label{eqpainleve}
	\gamma_{n}\Big(n+\lambda-z(\gamma_{n+1}+\gamma_{n})\Big)\Big(n-1+\lambda-z(\gamma_{n}+\gamma_{n-1})\Big)= \left(\dfrac{n}{2}-z\gamma_{n}\right)\left(\dfrac{n-1}{2}+\lambda-z\gamma_{n}\right).
\end{equation}
\end{teo}
\begin{proof}
Given the equation \eqref{eq6}, we deduce that $R_n=r_n$ and since $R_0=T_0=0$, equation \eqref{eq2} implies that $T_{n}=-t_n$. Additionally, taking into account equations \eqref{eq3} and \eqref{eq5}, we can conclude that
$$T_{n}=\dfrac{n}{2}-z\gamma_n,\quad\quad R_{n}=n+\lambda-z(\gamma_{n+1}+\gamma_n).$$
Now, from \eqref{eq1} and \eqref{eq5} for $j=0,1\ldots$
\begin{align*}
	\gamma_{j+1}R_{j+1}R_j-\gamma_{j}R_jR_{j-1}&=\left(\dfrac{1}{2}-z(\gamma_{j+1}-\gamma_j)\right)R_j\\
	&=(T_{j+1}-T_{j})R_j\\
	&=(T_{j+1}-T_{j})\left(T_{j+1}+T_{j}+\lambda-\dfrac{1}{2}\right)\\
	&=T_{j+1}^2-T_{j}^2+\left(\lambda -\dfrac{1}{2}\right)(T_{j+1}-T_j).
\end{align*}
A telescopic sum of the above equation yields
\begin{equation}\label{RTrela}\gamma_{n+1}R_{n+1}R_n=T_{n+1}\left(T_{n+1}+\lambda-\dfrac{1}{2}\right).\end{equation}
\end{proof}	
\begin{coro}\label{coro_simetricos_painleve}
If we define
\begin{equation}\label{gn}
	g_n=\dfrac{n}{2}+\dfrac{\lambda}{2}-\dfrac{1}{4}-z\gamma_n, \quad n\geq0, \end{equation} then \eqref{eqpainleve} becomes
\begin{equation}\label{painlevesimetricos}
	\left(\dfrac{n}{2}+\dfrac{\lambda}{2}-\dfrac{1}{4}-g_n\right)(g_{n+1}+g_n)(g_{n}+g_{n-1})=zg_n^2-z\left(\dfrac{\lambda}{2}-\dfrac{1}{4}\right)^2, \quad n\geq 1.
\end{equation}
Moreover, since $R_n$ is positive, $g_n+g_{n+1}\geq 0.$
Notice that it is a modified discrete Painlev\'e II equation, see \cite{Ramani}. When $\lambda=1/2$ we recover equation (56) in~\cite[Corollary~12]{DM22}.
\end{coro}
%\begin{remark}
%	The weight function $w(x,\lambda, z)=(1-x^2)^{(\lambda-1/2)} e^{-zx^2}$	belongs to Szeg\"o’s class, since 	the integral $$ \int^1_{
%	-1}\dfrac{\log w(x,\lambda,  z)}{\sqrt{1-x^2}}dx
%=-\dfrac{\pi}{2}\left[(2\lambda-1)\log(4)+z\right],$$
%and therefore, the following asymptotic property follows (see \cite{Walter1})
%\begin{equation*}%\label{asim}
% \lim_{n\to\infty}\gamma_{n}(z)=\dfrac{1}{4}.
%\end{equation*}
%Taking into account the above,  for $x$ on compact sets of $\mathbb{R}\setminus[-1,1]$ we get \cite{Walter2}
%$$\lim_{n\to\infty}P_n(x,z,\lambda)^{1/n}=\Phi(x),$$
%where $\Phi(x)$ is the solution of the quadratic equation
%$$ \Phi^2 - x\Phi + \dfrac{1}{4} = 0,\quad\longrightarrow \quad\Phi_{\pm}(x)=\dfrac{x\pm\sqrt{x^2-1}}{2}.$$
%From the above we deduce that
%$$
%P_n(x,z,\lambda) \sim \Phi_+(x)^n,\quad n\to\infty.$$

%Note that for $x \in (-1, 1)$
%$$\Phi_{\pm}(x)=\dfrac{x\pm i\sqrt{1-x^2}}{2}$$
%and therefore setting $x = \cos(\theta)$ we have
%$P_n\sim  \dfrac{1}{2^{n-1}}T_n(x)$ where $T_n(x)=\cos(n\theta)$ is the $n$-th Chebyshev polynomial of the first kind \cite{Chi}
%we see that
%$$P_n (x, z, \lambda)\sim \dfrac{1}{2^{n-1}}T_n(x)\quad n\to \infty$$
%\end{remark}
%%%%%%%%%%%%%%%%%%%%%%%%%%%%%%%%%%%%%%%%%%%%%	
\section{Holonomic Differential equation}\label{section5}
\begin{teo}[Structure relation] The  MOPS $(P_n)_{n\geq 0},$  orthogonal with respect to the linear functional $\sn,$  satisfies the differential-recurrence relation
\begin{equation}\label{estructureeq}
	\begin{aligned}
		\phi(x)\partial_xP_{n+1}(x)&=-(n+1)P_{n+2}(x)+b_{n+1}P_{n}(x)+a_{n}P_{n-2}(x),
	\end{aligned}
\end{equation}
where $\phi(x)=1-x^2$, $\gamma_0=b_0=a_0=a_1=0$ and
\begin{equation}\label{dttrr}
	\begin{aligned}
		b_{n+1}&=\left[n+1+2\lambda+2z-2z(\gamma_{n+2}+\gamma_{n+1}+\gamma_{n})\right]\gamma_{n+1},\quad n\geq0\\
		a_n&=-2z\gamma_{n+1}\gamma_{n}\gamma_{n-1}, \quad n\geq1.
	\end{aligned}
\end{equation}
\end{teo}
\begin{proof}
Taking into account that $(P_k)_{k=0}^{n+2}$ is a basis for the linear sequence of polynomials of degree less than or equal to $n+2$, then
$$\phi(x)\partial_xP_{n+1}(x)=-(n+1)P_{n+2}(x)+\sum_{k=0}^{n+1}\alpha_{n+1,k}P_k(x).$$
Using the orthogonality property and Pearson equation \eqref{pearson_eq}
\begin{equation}\label{def.alpha}
	\begin{aligned}
		\alpha_{n+1,k}h_k&=\prodint{\sn,\phi(\partial_xP_{n+1})P_k}=\prodint{\sn,\phi\partial_x(P_{n+1}P_k)}-\prodint{\sn,\phi P_{n+1}(\partial_xP_k)}\\
		&=\prodint{\sn,\psi P_{n+1}P_k}-\prodint{\sn,\phi  P_{n+1}(\partial_xP_k)}.
	\end{aligned}
\end{equation}
This implies that $\alpha_{n+1,k}=0$ for all $0\leq k<n-2.$ Moreover, since $\psi(x)=(1+2\lambda+2z)x-2zx^3$ is an odd function we get
$$\phi(x)\partial_xP_{n+1}(x)=-(n+1)P_{n+2}(x)+\alpha_{n+1,n}P_{n}(x)+\alpha_{n+1,n-2}P_{n-2}(x).$$
Now, from \eqref{eqq2} and \eqref{eqq3} we get
\begin{equation*}
	\begin{aligned}
		\alpha_{n+1,n}h_{n}&=\prodint{\sn,\psi P_{n+1}P_{n}}-\prodint{\sn,\phi  P_{n+1}(\partial_xP_{n})}\\&=\left[n+1+2\lambda+2z-2z(\gamma_{n+2}+\gamma_{n+1}+\gamma_{n})\right]\gamma_{n+1}h_{n},\\
		\alpha_{n+1,n-2}h_{n-2}&=\prodint{\sn,\psi P_{n+1}P_{n-2}}-\prodint{\sn,\phi  P_{n+1}P^\prime_{n-2}}=-2zh_{n+1}.
	\end{aligned}
\end{equation*}
Taking into account that $\gamma_n=\dfrac{h_n}{h_{n-1}}$, then
$$\alpha_{n+1,n-2}=-2z\dfrac{h_{n+1}}{h_{n-2}}=-2z\dfrac{h_{n+1}}{h_{n}}\dfrac{h_{n}}{h_{n-1}}\dfrac{h_{n-1}}{h_{n-2}}=-2z\gamma_{n+1}\gamma_{n}\gamma_{n-1},$$
and the result follows.
\end{proof}
Next, we deal with ladder operators associated with our MOPS (see \cite{MR1616931,Chen05}).
\begin{pro}\label{proL}
For $n\in\N$ let $L_n$ be the operator defined by
\begin{equation}\label{operatorL}
	L_n=A_n(x;z)\partial_x-B_{n}(x;z),
\end{equation}
where
\begin{equation*}
	\begin{aligned}
		A_n(x;z)=\dfrac{\phi(x)}{C_n(x;z;\lambda)},& &B_n(x;z)=\dfrac{\delta_n(x;z)}{C_{n}(x;z;\lambda)},
	\end{aligned}
\end{equation*}
with
$C_{n}(x;z;\lambda)=2\gamma_{n+1}[n+1+\lambda+z-z(\gamma_{n+2}+\gamma_{n+1})-zx^2]$ and $\delta_n(x;z)=(2z\gamma_{n+1}-n-1)x.$  Then
$$L_nP_{n+1}=P_n,\qquad n\geq0.$$
\end{pro}
\begin{proof}
Using the structure relation \eqref{estructureeq} and \eqref{ttrrs} we have
\begin{multline*}
	\phi (x)\partial_xP_{n+1}(x)=\\-(n+1)(xP_{n+1}(x)-\gamma_{n+1}P_n(x))+\alpha_{n+1,n}P_n(x)+\alpha_{n+1,n-2}P_{n-2}(x).
\end{multline*}
Notice also that
\begin{multline}\gamma_n\gamma_{n-1}P_{n-2}(x)=\gamma_n(xP_{n-1}(x)-P_{n}(x))=x^2P_n (x)-xP_{n+1}(x)-\gamma_{n}P_{n}(x).\end{multline}
So, remembering the values of $\alpha_{n+1,n+1}$ and $\alpha_{n+1,n}$ in \eqref{def.alpha} and taking into account the above, we get
\begin{multline*}
	\left[\phi(x)\partial_x+(n+1-2z\gamma_{n+1})x\right]P_{n+1}(x)=
	\\\hfill=[(n+1)\gamma_{n+1}+\alpha_{n+1,n}-2z\gamma_{n+1}x^2+2z\gamma_{n}\gamma_{n+1}]P_n(x)
	\\\hfill=2[n+1+\lambda+z-z(\gamma_{n+2}+\gamma_{n+1})-zx^2]\gamma_{n+1}P_n(x)
\end{multline*}
and the result follows.
\end{proof}
As a direct consequence, we can find a second order linear differential equation that the SMOP $\{P_{n}\}_{n\geq0}$ satisfies.
\begin{pro}%\label{DnP} 	
For a nonnegative integer $n\geq 0$, let $ D_{n+1}$ be the second order linear differential operator
\begin{equation}\label{Dn}\begin{aligned}
		&D_{n+1}=\\
		&C_{n}(x)\phi^2(x)\partial_{xx}-2\phi(x)\left[2zx\gamma_{n+1}(1-x^2)+x\left( \lambda+z+1/2-z{x}^{2}  \right)C_{n}(x)\right]\partial_{x}\\&+\left[\left(4\,{x}^{4}z+ 2\gamma_{n+1}{x}^{2}z- \left( n+2\,
		w+2\,z+1 \right) {x}^{2} +1\right)  \left( n+
		1-2\,z\gamma_{ n+1}  \right)
		\right]C_{n}(x)\\
		&+\dfrac{C^{2}_{n}(x) C_{n-1}(x)}{\gamma_n}+8\,\gamma_{n+1}  \left(n+1-2\,z\gamma_{n+1}
		\right)z{x}^{2}\phi(x) .
\end{aligned}\end{equation}
Then $D_{n+1}P_{n+1}=0.$
\end{pro}
\begin{proof}
From \eqref{operatorL} we get
\begin{equation}\label{edpL3}
	\begin{aligned}
		\phi (x)  P^\prime_{n+1}(x)&=\delta_{n}(x)P_{n+1}(x)+C_{n}(x)P_{n}(x)\\
		&=\delta_{n}(x)P_{n+1}(x)+\dfrac{C_{n}(x)}{\gamma_{n+1}}\left[xP_{n+1}(x)-P_{n+2}(x)\right]\\
		&=\left(\delta_{n}(x)+\dfrac{C_{n}(x) }{\gamma_{n+1}}x\right)P_{n+1}(x)-\dfrac{C_{n} (x)}{\gamma_{n+1}}P_{n+2}(x)
	\end{aligned}
\end{equation}
as well as
$$P_n(x)=\dfrac{\phi(x) }{C_{n} (x)}P^\prime_{n+1}(x)-\dfrac{\delta_{n}(x)}{C_{n}(x) }P_{n+1}(x).$$
Taking derivatives in the last equation
\begin{multline}\label{edpL}
	\phi(x) P^\prime_n(x)=\phi\Bigg[\left(\dfrac{\phi(x) }{C_{n}(x)}\right)^{\prime} P^{\prime}_{n+1}(x)+\dfrac{\phi(x)} {C_{n}(x)}P^{\prime\prime}_{n+1}(x)\\\hfill-\left(\dfrac{\delta_n(x)}{C_{n}(x)}\right)^{\prime}P_{n+1}(x)-\dfrac{\delta_ {n}(x)}{C_{n}(x)}P^{\prime}_{n+1}(x)\Bigg].	
\end{multline}
On the other hand, a shift $n+1\to n$ in \eqref{edpL3} yields
\begin{equation}\label{edpL1}
	\begin{aligned}
		\phi(x)  P^\prime_{n}(x)&=\left(\delta_{n-1}(x)+\dfrac{C_{n-1}(x)}{\gamma_{n}}x\right)P_{n}(x)-\dfrac{C_{n-1}(x)}{\gamma_{n}}P_{n+1}(x)\\
		&=\left(\delta_{n-1}(x)+\dfrac{C_{n-1}(x)}{\gamma_{n}}x\right)\left(\dfrac{\phi(x) }{C_{n}(x)}P^\prime_{n+1}(x)-\dfrac{\delta_{n}(x)}{C_{n}(x)}P_{n+1}(x)\right)\\&\hfill-\dfrac{C_{n-1}(x)}{\gamma_{n}}P_{n+1}(x).
	\end{aligned}
\end{equation}
Thus, from \eqref{edpL} and \eqref{edpL1}  we get the result.
\end{proof}
%%%%%%%%%%%%%%%%%%%%%%%%%%%%%%%%%%%%%%%%%%%%%%%%%%%%
\begin{remark}	
If $\Pn=(P_0,P_1,\cdots)^T$, then  the recurrence relation \eqref{ttrrs} can be written as
\begin{equation*}%\label{jacobi}
	x\Pn=J\Pn \quad \text{with}\quad J= \begin{pmatrix}
		0&1&0\\
		\gamma_1&0&1&0	\\
		0&\gamma_2&0&1&0\\
		&0&\gamma_3&0&1&0\\
		&&\ddots&\ddots&\ddots&\ddots&\ddots\\
	\end{pmatrix}.
\end{equation*}
Moreover, taking into account  \eqref{estructureeq} and \eqref{dttrr}, we also have
\begin{equation*}%\label{jacobidiff}
	\phi(x)\partial_{x}\Pn=L\Pn\quad \text{with}\quad L= \begin{pmatrix}
		0&0&0\\
		b_1&0&-1&0\\
		0&b_2&0&-2&0\\
		a_2&0&b_3&0&-3&0\\
		0&a_2&0&b_3&0&-4&0\\	
		&\ddots&\ddots&\ddots&\ddots&\ddots&\ddots&\ddots\\
	\end{pmatrix}.
\end{equation*}
From the above, we get  $$J\partial_x\Pn=x\partial_x\Pn+\Pn$$ as well as
\begin{align*}
	\phi(x)J\partial_x\Pn&=JL\Pn\\	
	\phi(x)J\partial_x\Pn&=\phi(x)(x\partial_x\Pn+\Pn)=LJ\Pn+(I-J^2)\Pn.
\end{align*}
Thus,  the compatibility condition of the Lax pair
$$[J, L] \Pn= JL\Pn-LJ\Pn=(I-J^2)\Pn$$
gives rise to two nontrivial equations:
\begin{align}
	\label{pain1}-(n-1)\gamma_n+b_{n+1}-b_{n}+n\gamma_{n+1}&=1-(\gamma_{n}+\gamma_{n+1}),&\quad n\geq 0,\\
	\notag \gamma_{n+1}b_n+a_{n+1}-a_n-b_{n+1}\gamma_n&=-\gamma_{n}\gamma_{n+1},& \quad n\geq 1.
\end{align}
Note that \eqref{pain1} is equivalent to the Laguerre-Freud equation \eqref{111}.
\end{remark}
%%%%%%%%%%%%%%%%%%%%%%%%%%%%%%%%%%%%%%%%%%
\section{Variable $z$}\label{Section6}
In this section we will study the behavior of the moments $(\sn_n(z,\lambda))_{n\geq 0}$ and the coefficients  $(\gamma_{n}(z,\lambda))_{n\geq 1}$ of \eqref{ttrrs}, as functions of the variable $z$.
 
Differentiating \eqref{innerproduct} with respect to $z$, we have
$$\partial_z\prodint{\sn,p}=\prodint{\sn,\partial_zp}-\prodint{\sn,x^2p}.$$
In particular, taking $p(x)=x^n$ we get
\begin{equation}\label{pzmoments}
\partial_z\sn_n=-\sn_{n+2},\quad n\geq 0.	
\end{equation}
Moreover, using \cite[eq. 13.3.15]{OLB10}, and \eqref{moments}
$$\dfrac{\sn_{2n+2}}{s_{2n}}=\partial_z\ln\left(\pFq{1}{1}\left(n+1/2;\lambda+n+1,-z\right)\right),\qquad n\geq 0.$$
Using \eqref{pzmoments}, we can obtain a  first order ODE (in $z$) for the Stieltjes function $\Su(t; z;\lambda)$.
\begin{teo}
Let $\Su_{\sn}(t; z,\lambda)$ be the Stieltjes function defined in Proposition \eqref{Stielt}. Then
\begin{equation}\label{Spz}
	\partial_z\Su_{\sn}(t; z,\lambda)=-t^2\Su_{\sn}(t;z;\lambda)+t\sn_0(z;\lambda).
\end{equation}
Moreover,
\begin{multline*}
	\Su_{\sn}(t;z;\lambda)=\dfrac{\sqrt{\pi}\,\Gamma(1/2+\lambda)}{\Gamma(\lambda+1)}\times\\\left[\dfrac{2(\lambda+1)t}{(2(\lambda+1)t^2-1)}\left(\pFq{1}{1}\left(1/2;\lambda+1;z\right)-e^{-zt^2}\right)+
	\dfrac{1}{t}e^{-zt^2}\pFq{2}{1}\left(1/2,1;\lambda+1,t^{-1}\right)\right].
\end{multline*}
\end{teo}
\begin{proof}
From the definition of  $\Su_{\sn}$ and \eqref{pzmoments}
$$\partial_{z}\Su_{\sn}=-\sum_{n=0}^{\infty}\dfrac{\sn_{n+2}}{t^{n+1}}=-t^2\left(\sum_{n=0}^\infty\dfrac{\sn_{n+2}}{t^{n+3}}+\dfrac{\sn_0}{t}-\dfrac{\sn_0}{t}\right)=-t^2\Su_{\sn}+t\sn_0.$$	
On the other hand, taking into account \eqref{Spz} and \eqref{initialcondi} we get
\begin{equation}\label{aux}
	e^{zt^2}\Su_{\sn}(t;z;\lambda)=t\int_0^ze^{x^2r}\sn_0(r,\lambda)dr+C(t,\lambda).\end{equation}
In particular, from \eqref{poh} and \eqref{moments},
$$\begin{aligned}
	C(t,\lambda)&=\Su_{\sn}(t;0;\lambda)=\dfrac{\Gamma(\lambda+1/2)}{t}\sum_{n=0}^\infty \dfrac{\Gamma(n+1/2)}{\Gamma(\lambda+n+1)}t^{-n}\\&=\dfrac{\sqrt{\pi}\,\Gamma(\lambda+1/2)}{\Gamma(\lambda+1)t}\sum_{n=0}^\infty\dfrac{(1/2)_n(1)_n}{(\lambda+1)_n}\dfrac{t^{-n}}{n!}\\
	&=\dfrac{\sqrt{\pi}\,\Gamma(\lambda+1/2)}{\Gamma(\lambda+1)t}\pFq{2}{1}\left(1/2,1;\lambda+1,t^{-1}\right).	
\end{aligned}$$
Since \cite[Section 13.3 (ii)]{OLB10}
$$\partial_z\left[e^{zt^2}\,\pFq{1}{1}\left(a;b,-z\right)\right]=\left(t^2-\dfrac{a}{b}\right)e^{zt^2}\,\pFq{1}{1}\left(a;b,-z\right),$$
then $$\begin{aligned}
	\int_0^ze^{t^2r}\sn_0(r,\lambda)&dr=\int_0^ze^{t^2r}\sqrt{\pi}\dfrac{\Gamma(\lambda+1/2)}{\Gamma(\lambda+1)}\times\pFq{1}{1}\left(1/2;\lambda+1,-r\right)dr\\
	&=\sqrt{\pi}\dfrac{\Gamma(\lambda+1/2)}{\Gamma(\lambda+1)}\dfrac{2(\lambda+1)}{2(\lambda+1)t^2-1}\left[e^{zt^2}\pFq{1}{1}\left(1/2;\lambda+1,-z\right)-1\right].
\end{aligned}$$
Changing the above in \eqref{aux} the result follows.
\end{proof}
\begin{teo}
The functions $h_n(z)$, $\gamma_{n}(z)$ satisfy the Toda-type equations
\begin{equation}\label{partialzh}
	\partial_z\ln(h_n)=-\gamma_{n+1}-\gamma_{n},	
\end{equation}
\begin{equation}\label{partialzgamma}
	\partial_z\ln(\gamma_n)=\gamma_{n-1}-\gamma_{n+1}.	
\end{equation}
\end{teo}
\begin{proof}
Notice that since $\partial_zP_n$ is a polynomial of degree less than or equal to $n-2$, then from \eqref{rrttx^2} we get
$$\partial_zh_n=2\prodint{\sn,P_n\,\partial_zP_n}-\prodint{\sn,x^2P_n^2}=-(\gamma_{n+1}+\gamma_n)h_{n}$$  and taking into account $\gamma_n=h_n/h_{n-1}$ we obtain \eqref{partialzh}.
To show \eqref{partialzgamma}, we use \eqref{partialzh}  together
$$\partial_z\gamma_n=\dfrac{\partial_zh_n}{h_{n-1}}-\dfrac{\partial_zh_{n-1}}{h_{n-1}}\dfrac{h_n}{h_{n-1}}=\dfrac{\partial_zh_n}{h_{n}}\gamma_n-\dfrac{\partial_zh_{n-1}}{h_{n-1}}\gamma_n=(\gamma_{n-1}-\gamma_{n+1})\gamma_n.$$	
\end{proof}
Combining the Laguerre-Freud equation \eqref{eqpainleve}  and the Toda-type
equation for $h_n (z)$, we obtain the following result
\begin{pro}
The functions $h_n(z)$ satisfy
\begin{multline*}
	\Big[(n+\lambda)h_n+z\partial_zh_n\Big]\Big[(n-1+\lambda)h_{n-1}+z\partial_zh_{n-1}\Big]\\=\left(\dfrac{n}{2}h_{n-1}-zh_{n}\right)\left(\dfrac{2\lambda+n-1}{2}h_{n-1}-zh_{n}\right),\quad n\geq 1.
\end{multline*}
\end{pro}
\begin{proof}
Denote $\partial_z h_n=h_n^\prime$. Replacing \eqref{partialzh} in \eqref{eqpainleve}, we get
\begin{equation*} \dfrac{1}{h^2_{n-1}}\Big((n+\lambda)h_n+zh^\prime_n\Big)\Big((n-1+\lambda)h_{n-1}+zh^\prime_{n-1}\Big)= \left(\dfrac{n}{2}-z\gamma_{n}\right)\left(\dfrac{n-1}{2}+\lambda-z\gamma_{n}\right),
\end{equation*}
and multiplying by $h^2_{n-1}$ on both sides of the above equation, we get the result.
\end{proof}
\begin{teo} The function $g_n(z)$ defined in \eqref{gn} satisfies
\begin{multline}\label{gx}
	\dfrac{2}{z}\left(2g_n-z\right)^2\left[\left(n+\lambda-\dfrac{1}{2}-2g_n\right)\left(4g_n^2-\left(\lambda-\dfrac{1}{2}\right)^2\right)+2z\left(g_n^\prime\right)^2\right]=\hfill\\\hfill
	\left[ 4g_n\left(3g_n-n-\lambda+\dfrac{1}{2}\right)-(\lambda-1/2)^2-2zg_{n}^{\prime\prime}-g_n^{\prime}\right]^2.
\end{multline}
\end{teo}
\begin{proof}
First of all, note that using the change of variable \eqref{gn}, we get $$g_{n+1}=g_{n-1}-\dfrac{g_n^\prime}{\gamma_n},\quad \text{(where $\partial_zg_n=g_n^\prime$)}$$  and, therefore, \eqref{painlevesimetricos} can be rewritten as
\begin{equation}\label{s1}
	\begin{aligned}
		\dfrac{\left(2g_n\right)^2-\left(\lambda-1/2\right)^2}{4\gamma_n}&=\left(g_{n-1}-\dfrac{g_n^\prime}{\gamma_n}+g_n\right)\left(g_{n-1}+g_n\right)\\
		&=\left(g_{n-1}+g_n-\dfrac{g_n^\prime}{2\gamma_n}\right)^2-\left(\dfrac{g_n^\prime}{2\gamma_n}\right)^2\\
		&=\left(g_{n-1}-\dfrac{g_n^\prime}{2\gamma_n}\right)^2+2\left(g_{n-1}-\dfrac{g_n^\prime}{2\gamma_n}\right)g_n+g_n^2 -\left(\dfrac{g_n^\prime}{2\gamma_n}\right)^2.
	\end{aligned}
\end{equation}

On the other hand, from  \eqref{eqpainleve} and \eqref{111} we get
$$\begin{aligned}
	z\gamma_{n-2}\gamma_{n-1}=&-\dfrac{T_{n-1}(T_{n-1}+\lambda-1/2)}{R_{n-1}}+(n-2+\lambda-z\gamma_{n-1})\gamma_{n-1},   \\
	z\gamma_{n+2}\gamma_{n+1}=&[z+\lambda+n+1-z\gamma_{n+1}]\gamma_{n+1}-[z+\lambda+n-1-z(\gamma_{n-1}+\gamma_{n})]\gamma_{n}-\dfrac{1}{2}.
\end{aligned}$$
Hence, using \eqref{partialzgamma}
\begin{equation*}
	\begin{aligned}
		z\dfrac{\gamma_{n}^{\prime\prime}}{\gamma_n}&=z\dfrac{\gamma_n^\prime(\gamma_{n-1}-\gamma_{n+1})}{\gamma_{n}}+z\dfrac{\gamma_n(\gamma_{n-1}-\gamma_{n+1})^\prime}{\gamma_{n}}\\
		&=z\left(\dfrac{\gamma_{n}^{\prime}}{\gamma_n}\right)^2+z\left(\gamma_{n+2}\gamma_{n+1}+\gamma_{n-2}\gamma_{n-1}\right)-z\gamma_n(\gamma_{n+1}+\gamma_{n-1})
		\\&=z\left(\dfrac{\gamma_{n}^{\prime}}{\gamma_n}\right)^2-\dfrac{T_{n-1}(T_{n-1}+\lambda-1/2)}{R_{n-1}}-\dfrac{1}{2} +\left(n-2+\lambda-z\gamma_{n-1}\right)\gamma_{n-1}\\&
		+\left(z+n+    1+\lambda-z\gamma_{n+1}\right)\gamma_{n+1}-\left[z+n+\lambda-1-z(\gamma_n-\gamma_{n+1})\right]\gamma_n\\
		&=-2z\gamma_{n+1}\gamma_{n-1}-\dfrac{T_{n-1}(T_{n-1}+\lambda-1/2)}{R_{n-1}}-\dfrac{1}{2} +\left(n-2+\lambda\right)\gamma_{n-1}\\&
		+\left(z+\lambda+n+1\right)\gamma_{n+1}-\left[z+\lambda+n-1-z(\gamma_n-\gamma_{n+1})\right]\gamma_n.
	\end{aligned}
\end{equation*}
From  \eqref{RTrela}
$$\begin{aligned}-\dfrac{T_{n-1}(T_{n-1}+\lambda-1/2)}{R_{n-1}}&=-R_{n-1}+2T_n+(\lambda-1/2)-\dfrac{T_{n}(T_{n}+\lambda-1/2)}{R_{n-1}}\\
	&=-R_{n-1}+2T_n+(\lambda-1/2)-\gamma_{n}R_n\\
	&=\dfrac{1}{2}+z\gamma_{n-1}-\left[z+n+\lambda-z(\gamma_{n+1}+\gamma_n)\right]\gamma_n.\end{aligned},$$
Thus, using   \eqref{partialzgamma} and the fact that $z\dfrac{\gamma_n^\prime}{\gamma_n}=-1-\dfrac{g_n^\prime}{\gamma_n}$, we get
\begin{multline*}
	\dfrac{z^2}{2}\dfrac{\gamma_{n}^{\prime\prime}}{\gamma_n}=
	-z^2\gamma_{n+1}\gamma_{n-1}+\dfrac{z\left(n+\lambda+z-2\right)}{2}\gamma_{n-1}\hfill\\\hfill+\dfrac{z\left(n+\lambda+z+1\right)}{2}
	\gamma_{n+1}-z\left(z+n+\lambda-\dfrac{1}{2}-z\gamma_n\right)\gamma_n\\
	=-\dfrac{z(n+\lambda+z+1)}{2}\dfrac{\gamma_n^\prime}{\gamma_n}
	+z\left[n-1+\lambda+z-\dfrac{1}{2}-z\gamma_{n-1}+1+z\dfrac{\gamma_n^\prime}{\gamma_n}\right]\gamma_{n-1}\hfill\\ \hfill-z\left(z+n+\lambda-\dfrac{1}{2}-z\gamma_n\right)\gamma_{n}\\
	=\dfrac{(n+\lambda+z+1)}{2}\left(1+\dfrac{g_n^\prime}{\gamma_n}\right)
	\hfill\\+\left(\dfrac{n-1}{2}+\dfrac{\lambda}{2}-\dfrac{1}{4}+g_{n-1}+z-\dfrac{g_n^\prime}{\gamma_n}\right)\left(\dfrac{n-1}{2}+\dfrac{\lambda}{2}-\dfrac{1}{4}-g_{n-1}\right)\hfill\\ \hfill+\left(\dfrac{n}{2}+\dfrac{\lambda}{2}-\dfrac{1}{4}+g_n+z\right)\left(g_n-\dfrac{n}{2}-\dfrac{\lambda}{2}+\dfrac{1}{4}\right)\\
	=\dfrac{(n+\lambda+z+1)}{2}\left(1+\dfrac{g_n^\prime}{\gamma_n}\right)+\left(\dfrac{n-1}{2}+\dfrac{\lambda}{2}-\dfrac{1}{4}+\dfrac{z}{2}-\dfrac{g_n^\prime}{2\gamma_n}\right)^2\hfill\\ \hfill -\left(g_{n-1}+\dfrac{z}{2}-\dfrac{g_n^\prime}{2\gamma_n}\right)^2
	-\left(\dfrac{n}{2}+\dfrac{\lambda}{2}-\dfrac{1}{4}+\dfrac{z}{2}\right)^2+\left(g_n+\dfrac{z}{2}\right)^2\\
	=\dfrac{(n+\lambda+z+1)}{2}\left(1+\dfrac{g_n^\prime}{\gamma_n}\right)-\dfrac{1}{2}\left(n+\lambda+z-1-\dfrac{g_n^\prime}{2\gamma_n}\right)\left(1+\dfrac{g_n^\prime}{\gamma_n}\right)\hfill\\\hfill
	-\left(g_{n-1}-\dfrac{g_n^\prime}{2\gamma_n}\right)^2-z\left(g_{n-1}-\dfrac{g_n^\prime}{2\gamma_n}\right)-\dfrac{z^2}{4}+\left(g_n+\dfrac{z}{2}\right)^2\\
	=\dfrac{1}{2}\left(2+\dfrac{g_n^\prime}{2\gamma_n}\right)\left(1+\dfrac{g_n^\prime}{\gamma_n}\right)-\left(g_{n-1}-\dfrac{g_n^\prime}{2\gamma_n}\right)^2-z\left(g_{n-1}-\dfrac{g_n^\prime}{2\gamma_n}\right)+\left(g_n+z\right)g_n.
\end{multline*}
Considering
$$\dfrac{z^2}{2}\dfrac{\gamma_{n}^{\prime\prime}}{\gamma_n}=1+\dfrac{g_n^\prime}{\gamma_n}-\dfrac{z}{2}\dfrac{g_n^{\prime\prime}}{\gamma_n},$$
then
\begin{equation}\label{s2}
	-\dfrac{z}{2\gamma_n}g_n^{\prime\prime}= \dfrac{g_n^\prime}{4\gamma_n}\left(1+\dfrac{g_n^\prime}{\gamma_n}\right)-\left(g_{n-1}-\dfrac{g_n^\prime}{2\gamma_n}\right)^2-z\left(g_{n-1}-\dfrac{g_n^\prime}{2\gamma_n}\right)+\left(g_n+z\right)g_n.
\end{equation}
Summing \eqref{s1} and \eqref{s2}
\begin{equation*}\begin{aligned}
		&\dfrac{\left(2g_n\right)^2-\left(\lambda-1/2\right)^2}{4\gamma_n}-\dfrac{z}{2\gamma_n}g_n^{\prime\prime}\\ &\qquad=\dfrac{g_n^\prime}{4\gamma_n}\left(1+\dfrac{g_n^\prime}{\gamma_n}\right)+(2g_n-z)\left(g_{n-1}-\dfrac{g_n^\prime}{2\gamma_n}\right)+\left(2g_n+z\right)g_n-\left(\dfrac{g_n^\prime}{2\gamma_n}\right)^2\\
		&\qquad=\dfrac{g_n^\prime}{4\gamma_n}+2zg_n+
		\left(g_{n-1}-\dfrac{g_n^\prime}{2\gamma_n}+g_n\right)(2g_n-z),
	\end{aligned}
\end{equation*}
or, equivalently,
\begin{equation*}
	\begin{aligned}
		g_{n-1}-\dfrac{g_n^\prime}{2\gamma_n}+g_n=&\dfrac{4g_n(g_n-2z\gamma_n)-\left(\lambda-\dfrac{1}{2}\right)^2-2zg_{n}^{\prime\prime}-g_n^{\prime}}{4\gamma_n\left(2g_n-z\right)}\\
		=&\dfrac{4g_n\left(3g_n-n-\lambda+\dfrac{1}{2}\right)-\left(\lambda-\dfrac{1}{2}\right)^2-2zg_{n}^{\prime\prime}-g_n^{\prime}}{4\gamma_n\left(2g_n-z\right)}.
	\end{aligned}
\end{equation*}
Replacing the above in \eqref{s1}  we get the result.
\end{proof}

\begin{remark}
Setting $g_n(z)=\nu_n\left(2\,\sqrt{ z}\right)$ in \eqref{gx} and taking $z=t^2/4$ we get

\begin{multline*}
	\left(\dfrac{2A\,\nu_n(t)}{t}-\dfrac{t}{A}\right)^2\left[\left(\nu^\prime(t)\right)^2-4\nu(t)^3 -\sigma_1\nu_n^2(t) - 2\beta_1\nu(t)-\varepsilon
	\right]=\hfill\\\hfill
	\left[\nu_{n}^{\prime\prime}(t)- 6\nu(t)^2 -\sigma_1\nu_n(t)- \beta_1\right]^2,
\end{multline*}
where
$$A=2,\quad \sigma_1=-2\left(n+\lambda-\dfrac{1}{2}\right),\quad \beta_1=-\left(\lambda-\dfrac{1}{2}\right)^2,\quad \varepsilon=\dfrac{\sigma_1\beta_1}{2}.$$
This is an example of the fourth Chazy's equation \cite[A.8]{Cosg06}, which is related to the Painlevé V equation.
\end{remark}

%%%%%%%%%%%%%%%%%%%%%%%%%%%%%%%%%%%%%%%%%%
\section{Zeros}\label{section7}
%%%%%%%%%%%%%%%%%%%%%%%%%%%%%%%%%%%%%%%%%%%%
It is well known that the zeros of orthogonal polynomials with respect to a positive definite linear functional are real, simple, and located in the interior of the convex hull of the support \cite{Chi,GMM21,Gabor}.
With this in mind, let $\{x_{n,k}(z;\lambda)\}_{k=1}^{n}$ be the zeros of $P_{n}(x, ;z; \lambda)$ in increasing order, that is,
\begin{equation}\label{zerosP_n}
P_{n}(x_{n,k}(z;\lambda);z;\lambda)=0
\end{equation}
with
$$x_{n,1}(z;\lambda)<x_{n,2}(z;\lambda)<\cdots<x_{n,n}(z;\lambda).$$
Next we will focus our attention on an electrostatic interpretation of the zeros of the polynomial $P_{n}(x ;z; \lambda)$ in terms of the energy associated with a logarithmic potential and next we will study the dynamics of them in terms of the parameter~$z.$
%%%%%%%%%%%%%%%%%%%%%%%%%%%%%%%%%%%%%%%%%%%%%
\subsection{Electrostatic interpretation}
%%%%%%%%%%%%%%%%%%%%%%%%%%%%%%%%%%%%%%%%%%%%%%%%%

Evaluating the operator $D_{n+1}$ defined in \eqref{Dn} at $x = x_{n+1,k}$, we get
\begin{align*}	\left.\dfrac{\partial_x^2P_{n+1}}{\partial_xP_{n+1}}\right|_{x=x_{n+1,k}}=\dfrac{4z\gamma_{n+1}x_{n+1,k}}{C_{n}(x_{n+1,k};z;\lambda)}+\dfrac{2x(\lambda+z+1/2-zx_{n+1,k}^2)}{\phi(x_{n+1,k})},
\end{align*}
where $C_{n}(x;z)$ was defined in Proposition \ref{proL}. The previous expression reads
\begin{multline}
\left.\dfrac{\partial_x^2P_{n+1}}{\partial_xP_{n+1}}\right|_{x=x_{n+1,k}}
=2zx_{n+1,k}-\dfrac{\lambda+1/2}{x_{n+1,k}-1}-\dfrac{\lambda+1/2}{x_{n+1,k}+1}\\-\dfrac{1}{x_{n+1,k}-\beta_{n+1}}-\dfrac{1}{x_{n+1,k}+\beta_{n+1}},
\end{multline}
where $$\beta^2_{n+1}=\dfrac{n+1+\lambda}{z}+1-(\gamma_{n+1}+\gamma_{n+2}).$$
Observe that the weight function associated with $\sn$ is $w(x)=(1-x^2)^{\lambda-1/2}e^{-zx^2}=e^{(\lambda-1/2)\ln(1-x^2)-zx^2}.$ If we define $v(x)=-(\lambda-1/2)\ln(1-x^2)+zx^2$ (the external potential \cite[Section~3.5]{Mourad}), then
\begin{multline}\label{equielect1}
\left.\dfrac{\partial_x^2P_{n+1}}{\partial_xP_{n+1}}\right|_{x=x_{n+1,k}}=-\dfrac{1}{x_{n+1,k}-\beta_{n+1}}-\dfrac{1}{x_{n+1,k}+\beta_{n+1}}\\
-\dfrac{1}{x_{n+1,k}-1}-\dfrac{1}{x_{n+1,k}+1}+v^\prime (x_{n+1,k}).
\end{multline}
\begin{teo}
The zeros of $P_{n+1}(x; z;\lambda)$ are the equilibrium points of $n+1$ unit charged particles located in the interval $(-1, 1 )$ under the influence of the  external potential
$$\ln\dfrac{1}{|x-\beta_{n+1}|}+\ln\dfrac{1}{|x+\beta_{n+1}|}
+\hat{v}(x), \quad -1<x<1,$$
where $\hat{v}(x)=- (\lambda + 1/2) \ln (1-x^2) + z x^2.$
\end{teo}
\begin{proof} If $P_{n+1}(x)=\prod_{k=1}^{n+1}(x-x_{n+1,k}),$
then \cite[Chapter 10]{GMM21}, \cite{Mourad} we have
\begin{equation*}%\label{cociente}
	\sum_{\substack{j=1 \\ j\ne k} }^{n+1}\dfrac{2}{x_{n+1,j}-x_{n+1,k}}=-\dfrac{\partial^2_xP_{n+1}(x_{n+1,k})}{\partial_xP_{n+1} (x_{n+1,k})}.
\end{equation*}
Taking into account  \eqref{equielect1}
\begin{multline*}
	\sum_{\substack{j=1 \\ j\ne k} }^{n+1}\dfrac{2}{x_{n+1,j}-x_{n+1,k}}-\dfrac{1}{x_{n+1,k}-\beta_{n+1}}-\dfrac{1}{x_{n+1,k}+\beta_{n+1}}\\
	-\dfrac{1}{x_{n+1,k}-1}-\dfrac{1}{x_{n+1,k}+1}+v^\prime (x_{n+1,k})=0.
\end{multline*}
It is not difficult to see that the above system of equations can be written as
$$\dfrac{\partial E}{\partial \,x_{k}}=0, \quad 1\le k \le n+1,$$
where
\begin{align*}
	E(&x_{1},\ldots, x_{n+1})=  2\sum_{1\le j < k \le n+1}\ln\dfrac{1}{|x_{j}-x_{k}|}+\\
	&\sum_{k=1}^{n+1}\left[\ln\dfrac{1}{|x_{k}-\beta_{n+1}|}+\ln\dfrac{1}{|x_{k}+\beta_{n+1}|}
	+\ln\dfrac{1}{|x_{k}-1|}+\ln\dfrac{1}{|x_{k}+1|}+v (x_{k})\right].
\end{align*}
\end{proof}
Notice that it involves the "natural" potential $v$ associated with the weight functions, positive unit charges located at the ends of the support of the weight, and two charges at the extra points $\pm\beta_{n}$. Since $$\beta_{n}=\sqrt{\dfrac{n+\lambda}{z}+1-(\gamma_{n}+\gamma_{n+1})}$$ and the notation in Corollary \ref{coro_simetricos_painleve} we get $z(\beta_{n}^2-1)=g_n+g_{n+1}=R_{n+1}\geq 0$. Therefore, $|\beta_{n}| \geq 1.$ Therefore, for all $n\geq1$  we conclude that $\beta_{n}$ is outside of the interval~$(-1.1)$.

\subsection{Dynamical behavior of zeros}
Now we focus our attention on the behavior of the zeros of $P_n(x,z)$ in terms of the variable $z$. In such a direction, we follow the approach given in \cite{Ismail-Wen11}.

First of all, notice that equation \eqref{partialzgamma} can be rewritten in terms of a Lax pair $(\Gamma_1(z),{\Gamma_1^2}_-(z)) $, that is,
\begin{align}\label{lax1}
\dot{\Gamma}_1(z)={\Gamma_1^2}_- {\Gamma_1} - {\Gamma_1} \, {\Gamma_1^2}_-,
\end{align}
where $\dot{q}(z):=\partial_zq(z)$ and ${\Gamma_1}(z)$ and ${\Gamma_1^2}_-(z)$ are matrices given by
\begin{align*}
\small
{\Gamma_1}(z)=
\begin{pmatrix}
	0 & 1%_{N}
	& & & \\
	\gamma_1(z) & 0 & 1%_{N}
	& & \\
	0 & \gamma_2(z) & 0 &1%_{N}
	\\
	& 0 & \gamma_3(z) & 0 & \ddots \\
	& & \ddots & \ddots & \ddots
\end{pmatrix}
, &&
\small
{\Gamma_1^2}_-(z)
=
\begin{pmatrix}
	0 &0 &\cdots\\
	0 &0 &\cdots\\
	\gamma_2(z)\gamma_1(z)&0 & \cdots\\
	0 &\gamma_3(z)\gamma_2(z)& \\
	\vdots &\vdots &\ddots \\
\end{pmatrix}.
\end{align*}
Notice also that the three-term recurrence relation \eqref{rr} can be written as
\begin{equation}\label{lax_z}
x\Pn=\Gamma_1\Pn,\qquad \Pn=(P_0,P_1,\cdots)^T.
\end{equation} Taking derivatives in \eqref{lax_z} with respect to $z$ and in view of \eqref{lax1} we get
\begin{align*}
(\Gamma_{1-}^2\Gamma_{1}-\Gamma_{1}\Gamma_{1-}^2)\Pn(x,z)
-(xI-{\Gamma_1})\dot{\Pn}(x,z)& =0,\\
x\Gamma_{1-}^2-\Gamma_{1}\Gamma_{1-}^2\Pn(x,z)
-(xI-{\Gamma_1})\dot{\Pn}(x,z)& =0,\\
(xI-\Gamma_1)(\Gamma_{1-}^2\Pn(x,z)-\dot{\Pn}(x,z))&=0.
\end{align*}
Since for each $z$, $( P_n(x; z))_{n\geq 0}$ is a basis of $\mathbb{P}$, we conclude $$\dot{\Pn}(x,z)=\Gamma_{1-}^2\Pn(x,z).$$
Hence,
\begin{pro}
Let  $ (P_n)_{n\geq 0}$ be the sequence of monic orthogonal polynomials with respect to $\sn$. Then, for every $n\geq 2$,
\begin{equation*}%\label{Pprime}
	\dot{P}_n(x,z)=\gamma_{n}\gamma_{n-1}P_{n-2}(x,z).
\end{equation*}
\end{pro}
If we take the derivative of \eqref{zerosP_n} with respect to $z,$ then we get 	$$\partial_xP_n(x,z)\left|_{x=x_{n,k}(z)}\right.\partial_z x_{n,k}(z)\,+\dot {P}_n(x,z)\left|_{x=x_{n,k}(z)}\right.=0,$$
or, equivalently,
\begin{equation}\label{diffy}
\dot{x}_{n,k}(z)=-\dfrac{\dot {P}_n(x,z)\left|_{x=x_{n,k}(z)}\right.}{\partial_x\,P_n(x,z)\left|_{x=x_{n,k}(z)}\right.}.
\end{equation}
Furthermore, from Proposition \ref{proL} we deduce that
\begin{equation}\label{diffP}
\partial_xP_{n}(x,z)\left |_{x=x_{n,k}(z)}\right.=\dfrac{P_{n-1}(x_{n,k}(z),z)}{A_{n-1}(x_{n,k}(z),z)},
\end{equation}
where $A_{n-1}(x,z)$ was defined in \eqref{operatorL}.
\begin{pro}
Let $(x_{n,k}(z))_{k=1}^n$ be the zeros of $P_{n}(x,z)$. Then they satisfy the differential equation \begin{equation*}%\label{diifx}
	\dot{x}_{n,k}(z)=-  \dfrac{x_{n,k}(z)}{2z} \dfrac{1- x^{2}_{n,k}(z)}{ \beta^{2}_{n} - x^{2}_{n,k} (z)}.
\end{equation*}
\end{pro}
\begin{proof}
Taking into account \eqref{rr}, \eqref{operatorL}, \eqref{diffy} and \eqref{diffP}
$$\begin{aligned} \dot{x}_{n,k}(z)&=-A_{n-1}(x_{n,k}(z),z)\,\dfrac{\dot {P}_n(x,z)\left|_{x=x_{n,k}(z)}\right.}{P_{n-1}(x_{n,k}(z),z)}\\
	&=-A_{n-1}(x_{n,k}(z),z)\,\dfrac{\gamma_n\gamma_{n-1}P_{n-2}(x_{n,k}(z),z)}{P_{n-1}(x_{n,k}(z),z)}\\
	&=-A_{n-1}(x_{n,k}(z),z)\,\dfrac{\gamma_{n}x_{n,k}(z)P_{n-1}(x_{n,k}(z),z)}{P_{n-1}(x_{n,k}(z),z)},\end{aligned}$$

and the result follows from the expression of $A_{n-1} (x;z)$  in \eqref{operatorL}.
\end{proof}
\begin{remark}
Notice that if  we denote $y_{n,k}(z)= x^{2}_{n,k}(z),$ then the above expression reads
$$\dot{y}_{n,k}\left( \dfrac{\beta^{2}_{n}}{ y_{n,k}(z)} + \dfrac{ \beta^{2}_{n}-1}{ 1- y_{n,k}(z)}\right)= - \dfrac{1}{z}.$$
As a straightforward consequence we get that  $\dot{y}_{n,k} <0,$ i. e.,  they are decreasing functions in terms of $z.$
\end{remark}
%%%%%%%%%%%%%%%%%%%%%%%%%%%%%%%%%%%%%%%%
%%%%%%%%%%%%%%%%%%%%%%%%%%%%%%%%%%%%%


\begin{thebibliography}{30}

\bibitem{Belmehdi} S. Belmehdi,\textit{ On semi--classical linear functionals of class s=1. Classification and integral representations.}  Indag. Math. (N.S.) \textbf{3} (1992), no. 3, 253--275.

\bibitem{Said} S. Belmehdi,  A. Ronveaux, \textit{Laguerre--Freud's equations for the recurrence coefficients of semi--classical orthogonal polynomials}. J. Approx. Theory \textbf{76}  (1994), no. 3,  351--368.

\bibitem{MR1083352} C. Brezinski, \textit{History of Continued Fractions and Padé Approximants.} Springer Series in Computational Mathematics, {\bf 12},  Springer--Verlag, Berlin, 1991.

		
\bibitem{MR1616931} Y. Chen M. E. H. Ismail, \textit{ Ladder operators and differential equations for orthogonal
polynomials.} J. Phys. A {\bf 30} (1997), no. 22, 7817--7829.
		
		
\bibitem{Chen05} Y. Chen, G. Pruessner, \textit{Orthogonal polynomials with discontinuous weights}. J. Phys. A\textbf{ 38} (2005), no. 12,  L191--L198.
		
%\bibitem{MR2676782}	Y. Chen, L. Zhang, \textit{Painlev\'{e} {VI} and the unitary Jacobi ensembles.} Stud. Appl. Math. {\bf125} (2010), no. 1, 91--112.
		
%\bibitem{Ted} T. S. Chihara, \textit{Generalized Hermite polynomials}, Doctoral Dissertation. Purdue University. 1955.

\bibitem{Cosg06} C. M. Cosgrove, \textit{Chazy’s second-degree Painlevé equations,} J. Phys. A {\bf 39} (2006), 11955-11971.
		
\bibitem{Chi} T. S. Chihara, \textit{An Introduction to Orthogonal Polynomials}. In : Mathematics and its Applications Series, Vol. \textbf{13}. Gordon and Breach Science Publishers, New York--London--Paris, 1978.
		
%\bibitem{Deift99} P. A. Deift, \textit{Orthogonal polynomials and random matrices: a Riemann--Hilbert approach.}Courant Lecture Notes in Mathematics, {\bf 3}. New York University, Courant Institute of Mathematical Sciences, New York, American Mathematical Society, Providence, RI, 1999. viii+273 pp.
		
%\bibitem{MR4093808} D. Dominici, \textit{Power series expansion of a Hankel determinant}. Linear Algebra Appl. {\bf 601} (2020), 17--54.
		
%\bibitem{MR4136730} D. Dominici, \textit{Recurrence coefficients of Toda--type orthogonal polynomials I. Asymptotic analysis.} Bull. Math. Sci. {\bf10}, no.2, 2050003, 32 pp.

\bibitem{DGM23} D. Dominici, J. C. García--Ardila, F. Marcellán, \textit{Symmetrization process and truncated orthogonal polynomials.} arXiv:2307.09581v2 [math.CA]		
		
\bibitem{DM22} D. Dominici, F. Marcellán \textit{Truncated Hermite polynomials.} J. Difference Equ. Appl.  \textbf{29} (2023) no. 7, 701-732.

 \bibitem{Du76} M. Dupuis, J. Rys, H. F.  King,  \textit{Evaluation of molecular integrals over Gaussian basis functions.} J. Chem. Phys. \textbf{65} (1976) , 111-116.
		
%\bibitem{MR1885665} P. J. Forrester, N. S. Witte, \textit{Application of the {$\tau$}--function theory of {P}ainlev\'{e}
%equations to random matrices: {$\rm P_V$}, {$\rm P_{III}$}} Comm. Pure Appl. Math. {\bf 55} (2002), no. 6, 679--727.

\bibitem{GM23} J. C. García--Ardila, F. Marcellán,  \textit{A  note on Laguerre truncated polynomials and quadrature formula.}  	arXiv:2401.00752 [math.NA], 2024.

\bibitem{GMM21} J. C. García--Ardila, F. Marcellán, M. E. Marriaga, \textit{Orthogonal Polynomials and Linear Functionals. An Algebraic Approach and Applications}. European Mathematical Society (EMS), Zürich, 2021.

\bibitem{Ga82} W. Gautschi,  \emph{On generating orthogonal polynomials}. SIAM J. Sci. Statist. Comput. \textbf{3 }(1982),  282–317.

\bibitem{Ga04} W. Gautschi,  \emph{Orthogonal Polynomials: Computation and Approximation}. Clarendon Press, Oxford. 2004.

\bibitem{MR1164865} P. Henrici, \textit{Applied and Computational Complex Analysis Vol. \textbf{2}.} Special functions--integral transforms--asymptotics--continued fractions. Wiley Classics Lib. John Wiley \& Sons, Inc., New York, 1991.

\bibitem{Mourad} M. E. H. Ismail, \textit{Classical and Quantum Orthogonal Polynomials in One Variable}.  Encyclopedia of Mathematics and its Applications, \textbf{98}, Cambridge University Press, Cambridge, 2005.

\bibitem{Ismail-Wen11} M. E. H. Ismail,   W.-X. Ma, \textit{Equations of motion for zeros of orthogonal polynomials related to the Toda lattices.} Arab J. Math. Sci. {\bf 17} (2011), no. 1, 1--10.

\bibitem{K16} H. F. King, \emph{Strategies for evaluation of Rys roots and weights.} J. Phys. Chem. A \textbf{120} (2016), 9348–9351.
		
\bibitem{MR2442472} S. Khrushchev, \textit{Orthogonal Polynomials and Continued Fractions. From Euler's point of view.} Encyclopedia Mathematics and its Applications,  {\bf 122}, Cambridge University Press, Cambridge, 2008.

\bibitem{Lyu1} S. Lyu,  Y. Chen, \textit{Gaussian unitary ensembles with two jump iscontinuities, PDEs, and the coupled Painlev\'e II and IV systems.} Stud. Appl. Math. \textbf{146}(1) (2021) 118--138.
		
%\bibitem{Lyu2} S. Lyu, Y. Chen, S--X. Xu, \textit{Laguerre Unitary Ensembles with Jump Discontinuities, PDEs and the Coupled Painlev\'e V System}. arXiv:2202.00943v1[nlin.SI].
		
\bibitem{Magnus1} A. P. Magnus, \textit{ Painlev\'e--type  differential equations for the recurrence coefficients of semi--classical orthogonal polynomials.} In \textit{ Proceedings of the Fourth International Symposium on Orthogonal Polynomials and their Applications (Evian--Les--Bains, 1992)}. J. Comput. Appl. Math. {\bf 57} (1995), no. 1--2, 215--237.
		
\bibitem{Magnus2} A. P. Magnus,  \textit{Freud's equations for orthogonal polynomials as discrete Painlev\'e equations.} In \textit{Symmetries and integrability of difference equations (Canterbury, 1996)}, vol. \textbf{255},  London Math. Soc. Lecture Note Ser.  228--243. Cambridge Univ. Press, Cambridge, 1999.
		
		
%\bibitem{Ma85} P. Maroni, \textit{Une caract\'erisation des polyn\^omes orthogonaux semi--classiques}, C. R. Acad. Sci. Paris Sér. I Math. \textbf{301} (1985), no. 6, 269--272.
		
\bibitem{Ma87} P. Maroni, \textit{  Prol\'egom\`enes \`a  l'\'etude des polyn\^{o}mes orthogonaux semi--classiques.}.  Ann. Mat. Pura Appl. (4)\textbf{ 149 }(1987), 165--184.
		
\bibitem{Ma91} P. Maroni,  \textit{Une th\'eorie alg\'ebrique des polyn\^{o}mes orthogonaux. Application aux polyn\^{o}mes orthogonaux semi--classiques.}  In  \textit{Orthogonal polynomials and their applications (Erice, 1990)}, 95--130, IMACS Ann. Comput. Appl. Math., \textbf{9}, Baltzer, Basel, 1991.

\bibitem{Gradimir14} G. V. Milovanović, \emph{Orthogonal polynomials on the real line}. Chapter 11 In: \emph{Walter Gautschi: Selected Works with Commentaries}, Volume \textbf{2 }(C. Brezinski, A. Sameh, eds.),  3-16, Birkh\"{a}user, Basel. 2014.
		
\bibitem{Gradimir18} G. V. Milovanović, \textit{An efficient computation of parameters in the Rys quadrature formula}. Bull. Cl. Sci. Math. Nat. Sci. Math. \textbf{43} (2018), 39--64.

%\bibitem{MMR94}  G. V. Milovanović, D. S.  Mitrinović,  Th. M. Rassias, \textit{Topics in polynomials: extremal problems, inequalities, zeros.}
%World Scientific Publishing Co., Inc., River Edge, NJ, 1994.
		
\bibitem{Gradimir22}  G. V. Milovanović, N. Vasović,\textit{ Orthogonal polynomials and generalized Gauss--Rys quadrature formulae}. Kuwait J. Sci. \textbf{49} (2022), no. 1, 17 pp.
		
%\bibitem{MR914314} K. Okamoto, \textit{ Studies on the Painlev\'{e} equations. {II}. {F}ifth {P}ainlev\'{e} equation {$P_{\rm V}$}.} Japan. J. Math. (N.S.) {\bf 13} (1987), no. 1, 47--76
		
\bibitem{OLB10} F. W. J. Olver, D. W. Lozier, R. F. Boisvert, C. W. Clark, editors, \textit{NIST Handbook of Mathematical Functions}. U.S. Department of Commerce, National Institute of Standards and Technology, Washington, DC,  Cambridge University Press, Cambridge, 2010.
		
\bibitem{Ramani} A. Ramani, B. Grammaticos, \textit{Miura transforms for discrete Painlevé equations}. J. Phys. A \textbf{25} (1992), L633--L637.

\bibitem{Shizgal} B. D. Shizgal, \textit{ A novel Rys quadrature algorithm for use in the calculation of electron repulsion integrals}. Comput. Theor. Chemistry \textbf{1074} (2015), 178--184.
		
\bibitem{Gabor} G. Szeg\H{o}, \textit{Orthogonal Polynomials}. Amer. Math. Soc. Colloq. Public. vol \textbf{23}, Amer. Math. Soc. Providence RI, 1975. Fourth Edition.

%\bibitem{Toda89} M. Toda, \textit{Theory of nonlinear lattices.} Translated from the Japanese by the author. Springer Series in Solid--State Sciences, {\bf 20}. Springer--Verlag, Berlin--New York, 1981. x+205 pp.
%\bibitem{Walter2} W. Van Assche, \textit{Asymptotics for Orthogonal Polynomials and Three-Term Recurrences.} In: Nevai, P. (eds) Orthogonal  Polynomials. NATO ASI Series, vol. {\bf 294.} Springer, Dordrecht. 1990.

%\bibitem{Walter1} W. Van Assche, \textit{ Asymptotics for orthogonal polynomial}, vol. {\bf 1265}
%of “Lecture Notes in Mathematics”. Springer-Verlag, Berlin (1987).

\bibitem{Walter} W. Van Assche, \textit{Orthogonal Polynomials and Painlev\'e Equations.} Australian Mathematical Society Lecture Series, \textbf{27}. Cambridge University Press, Cambridge, 2018.

\bibitem{MR0025596} H. S. Wall, \textit{Analytic Theory of Continued Fractions.}
D. Van Nostrand Co., Inc., New York, NY,  1948.

%\bibitem{Wu} X. Wu,  S. Xu, \textit{Gaussian unitary ensemble with jump discontinuities and the coupled Painlev\'e $II$ and $ IV$ systems}, Nonlinearity \textbf{34 }(2021), 2070--2115.

\end{thebibliography}
\end{document}